\documentclass[10pt]{article}
\usepackage{cite}
\usepackage{mathrsfs}
\usepackage{amsfonts}
\usepackage{amsmath}
\usepackage{amsfonts,amssymb}
\usepackage{dsfont}
\usepackage{curves}
\usepackage{mathrsfs}
\usepackage{pifont}
\usepackage{amssymb}
\allowdisplaybreaks

\numberwithin{equation}{section}

\date{}

\def\BigRoman{\uppercase\expandafter{\romannumeral\number\count 255 }}
\def\Romannumeral{\afterassignment\BigRoman\count255=}

\begin{document}
\title{A spectral condition for spanning trees with restricted degrees in bipartite graphs
}
\author{\small  Jiancheng Wu$^{1}$, Sizhong Zhou$^{1}$\footnote{Corresponding author. E-mail address: zsz\_cumt@163.com (S. Zhou)},
Hongxia Liu$^{2}$\\
\small  $1$. School of Science, Jiangsu University of Science and Technology,\\
\small  Zhenjiang, Jiangsu 212100, China\\
\small  $2$. School of Mathematics and Information Science, Yantai University,\\
\small  Yantai, Shandong 264005, China\\
}

\maketitle
\begin{abstract}
\noindent Let $G$ be a graph and $T$ be a spanning tree of $G$. We use $Q(G)=D(G)+A(G)$ to denote the signless Laplacian matrix of $G$, where
$D(G)$ is the diagonal degree matrix of $G$ and $A(G)$ is the adjacency matrix of $G$. The signless Laplacian spectral radius of $G$ is denoted
by $q(G)$. A necessary and sufficient condition for a connected bipartite graph $G$ with bipartition $(A,B)$ to have a spanning tree $T$ with
$d_T(v)\geq k$ for any $v\in A$ was independently obtained by Frank and Gy\'arf\'as (A. Frank, E. Gy\'arf\'as, How to orient the edges of a graph?,
Colloq. Math. Soc. Janos Bolyai 18 (1976) 353--364), Kaneko and Yoshimoto (A. Kaneko, K. Yoshimoto, On spanning trees with restricted degrees,
Inform. Process. Lett. 73 (2000) 163--165). Based on the above result, we establish a lower bound on the signless Laplacian spectral radius $q(G)$
of a connected bipartite graph $G$ with bipartition $(A,B)$, in which the bound guarantees that $G$ has a spanning tree $T$ with $d_T(v)\geq k$ for
any $v\in A$.
\\
\begin{flushleft}

(2020) Mathematics Subject Classification: 05C50, 05C05, 05C70
\end{flushleft}
\end{abstract}

\section{Introduction}

Throughout this paper, we only discuss simple, undirected and connected graphs. Let $G=(V(G),E(G))$ be a graph, where $V(G)$ denotes its vertex set and
$E(G)$ denotes its edge set. For a vertex $v\in V(G)$, the neighborhood $N_G(v)$ of $v$ in $G$ is defined by $\{u\in V(G):uv\in E(G)\}$ and the number
$d_G(v)=|N_G(v)|$ is the degree of $v$ in $G$. For a vertex subset $S\subseteq V(G)$, we write $N_G(S)=\bigcup_{v\in S}N_G(v)$.

For a given graph $G$ with vertex set $V(G)=\{v_1,v_2,\cdots,v_n\}$, the adjacency matrix of $G$ is defined by $A(G)=(a_{ij})$, where $a_{ij}=1$
if two vertices $v_i$ and $v_j$ are adjacent in $G$, and $a_{ij}=0$ otherwise. Let $Q(G)=D(G)+A(G)$ denote the signless Laplacian matrix of $G$, where
$D(G)=\mbox{diag}\{d_G(v_1),d_G(v_2),\ldots,d_G(v_n)\}$ is the diagonal degree matrix of $G$. Let $\rho_1(G)\geq\rho_2(G)\geq\cdots\geq\rho_n(G)$ and
$q_1(G)\geq q_2(G)\geq\cdots\geq q_n(G)$ be the eigenvalues of $A(G)$ and $Q(G)$, respectively. In particular, the largest eigenvalue $\rho_1(G)$ of
$A(G)$ is called the adjacency spectral radius of $G$ and denoted by $\rho(G)$, and the largest eigenvalue $q_1(G)$ of $Q(G)$ is called the signless
Laplacian spectral radius of $G$ and denoted by $q(G)$. Some properties on spectral radius can be found in \cite{BOT,Oa,OPPZ,ZZS,ZZB,ZLt,ZSL1,ZZL,ZZ}.

Let $a$ and $b$ be two integers with $0\leq a\leq b$. Then a spanning subgraph $F$ of $G$ is called an $[a,b]$-factor of $G$ if $d_F(v)\in[a,b]$ for
any $v\in V(G)$. A spanning tree $T$ of a connected graph $G$ is called a spanning $k$-tree of $G$ if $d_T(v)\leq k$ for each $v\in V(G)$, that is,
the maximum degree of a spanning $k$-tree of $G$ is at most $k$, where $k\geq2$ is an integer. Obviously, a spanning $k$-tree of $G$ is also a connected
$[1,k]$-factor of $G$. A spanning tree having at most $k$ leaves is called a spanning $k$-ended tree. Kaneko \cite{Ks} introduced the concept of leaf
degree of a spanning tree. Let $T$ denote a spanning tree of a connected graph $G$. The number of leaves adjacent to a vertex $v$ in $T$ is called the
leaf degree of $v$. Furthermore, the maximum leaf degree among all the vertices in $T$ is called the leaf degree of $T$. The minimum of distances between
any two leaves in $T$ is called the leaf distance of $T$.

Lots of scholars investigated the existence of spanning trees under some specified conditions. Ota and Sugiyama \cite{OS} posed a sufficient condition
for a graph to contain a spanning $k$-tree via the condition on forbidden subgraphs. Kyaw \cite{Kyaw} obtained a degree and neighborhood condition for
the existence of a spanning $k$-tree in a graph. Win \cite{Win} showed some results on the existence of a spanning $k$-tree in a graph. Zhou and Wu
\cite{ZW} provided an upper bound on the distance spectral radius in a graph $G$ to ensure the existence of a spanning $k$-tree in $G$. Zhou, Zhang and
Liu \cite{ZZL1} studied the relation between the spanning $k$-tree and the distance signless Laplacian spectral radius in a graph and claimed an
upper bound on the distance signless Laplacian spectral radius in a graph $G$ to ensure the existence of a spanning $k$-tree in $G$. Broersma and Tuinstra
\cite{BT} presented a degree condition for a graph to contain a spanning $k$-end tree. Ao, Liu and Yuan \cite{ALY} obtained tight spectral conditions
to guarantee a graph to have a spanning $k$-end tree, and also posed tight spectral conditions for the existence of a spanning tree with leaf degree at
most $k$ in a graph. Zhou, Sun and Liu \cite{ZSL2} showed the upper bounds for the distance spectral radius (resp. the distance signless Laplacian
spectral radius) of a graph $G$ to guarantee that $G$ has a spanning tree with leaf degree at most $k$. Wu \cite{Wc} gave a lower bound on the size of
a graph $G$ to guarantee that $G$ has a spanning tree with leaf degree at most $k$, and established a lower bound on the spectral radius of a graph $G$
to ensure that $G$ has a spanning tree with leaf degree at most $k$. Kaneko, Kano and Suzuki \cite{KKS} posed a sufficient condition for a graph to
have a spanning tree with leaf distance at least 4. Erbes, Molla, Mousley and Santana \cite{EMMS} investigated the existence of spanning trees with
leaf distance at least $d$, where $d\geq4$ is an integer. Wang and Zhang \cite{WZ} showed an $A_{\alpha}$-spectral radius condition for the existence of
a spanning tree with leaf distance at least 4 in a graph. For more results on spanning subgraphs, we refer the reader to \cite{KL,KOPR,GWW,ZPX1,ZXS,Za1}.

Let $G$ be a bipartite graph with bipartition $(A,B)$. Let $K_{m,n}$ denote the complete bipartite graph with bipartition $(A,B)$, where $|A|=m$ and
$|B|=n$. Given two bipartite graphs $G_1=(A_1,B_1)$ and $G_2=(A_2,B_2)$, let $G_1\nabla G_2$ denote the graph obtained from $G_1\cup G_2$ by adding all
possible edges between $A_2$ and $B_1$.

For bipartite graphs, Kano, Matsuda, Tsugaki and Yan \cite{KMTY} provided a degree condition for a connected bipartite graph to contain a spanning
$k$-ended tree. Frank and Gy\'arf\'as \cite{FG}, Kaneko and Yoshimoto \cite{KY} independently studied the existence of a spanning tree $T$ with
$d_T(v)\geq k$ for any $v\in A$ in a connected bipartite graph $G$ with bipartition $(A,B)$ and obtained a necessary and sufficient condition for
the connected bipartite graph $G$ to have a spanning tree $T$ with $d_T(v)\geq k$ for any $v\in A$. Motivated by \cite{FG,KY} directly, it is natural
and interesting to put forward some sufficient conditions to guarantee that a connected bipartite graph with bipartition $(A,B)$ has a spanning tree
$T$ with $d_T(v)\geq k$ for any $v\in A$ with respect to the spectral radius. Our main result is shown as follows.

\medskip

\noindent{\textbf{Theorem 1.1.}} Let $k$, $m$ and $n$ be three integers with $k\geq3$, $m\geq3$ and $n\geq(k-1)m+1$, and let $G$ be a connected bipartite
graph with bipartition $A\cup B$, where $|A|=m$ and $|B|=n$. If
$$
q(G)\geq q(K_{1,k-1}\nabla K_{m-1,n-k+1}),
$$
then $G$ contains a spanning tree $T$ with $d_T(v)\geq k$ for any $v\in A$, unless $G=K_{1,k-1}\nabla K_{m-1,n-k+1}$.

\medskip

\section{Some preliminaries}

In this section, we introduce some lemmas, which will be used in the proofs of our main results.

\medskip

\noindent{\textbf{Lemma 2.1}} (Shen, You, Zhang and Li \cite{SYZL}). Let $G$ be a connected graph, and let $H$ be a subgraph of $G$. Then
$$
q(H)\leq q(G),
$$
where the equality holds if and only if $H=G$.

\medskip

Let $M$ be a real symmetric matrix whose rows and columns are indexed by $V=\{1,2,\cdots,n\}$. Assume that $M$, with respect to the partition
$\pi:V=V_1\cup V_2\cup\cdots\cup V_m$, can be written as
\begin{align*}
M=\left(
  \begin{array}{cccc}
    M_{11} & M_{12} & \cdots & M_{1m}\\
    M_{21} & M_{22} & \cdots & M_{2m}\\
    \vdots & \vdots & \ddots & \vdots\\
    M_{m1} & M_{m2} & \cdots & M_{mm}\\
  \end{array}
\right),
\end{align*}
where $M_{ij}$ denotes the submatrix (block) of $M$ formed by rows in $V_i$ and columns in $V_j$. Let $q_{ij}$ denote the average row sum of
$M_{ij}$. Then matrix $M_{\pi}=(q_{ij})$ is called the quotient matrix of $M$. If the row sum of each block $M_{ij}$ is a constant, then the
partition is equitable.

\medskip

\noindent{\textbf{Lemma 2.2}} (You, Yang, So and Xi \cite{YYSX}).  Let $M$ be a real symmetric matrix with an equitable partition $\pi$, and
let $M_{\pi}$ be the corresponding quotient matrix. Then every eigenvalue of $M_{\pi}$ is an eigenvalue of $M$. Furthermore, if $M$ is
nonnegative, then the largest eigenvalues of $M$ and $M_{\pi}$ are equal.

\medskip

Frank and Gy\'arf\'as \cite{FG}, and Kaneko and Yoshimoto \cite{KY} put forward a necessary and sufficient condition for bipartite graphs to have
spanning trees with restricted degrees, independently.

\medskip

\noindent{\textbf{Lemma 2.3}} (Frank and Gy\'arf\'as \cite{FG}, Kaneko and Yoshimoto \cite{KY}). Let $G$ be a connected bipartite simple graph with
bipartition $A\cup B$, and $f:A\longrightarrow\{2,3,4,\ldots\}$ be a function. Then $G$ contains a spanning tree $T$ such that $d_T(v)\geq f(v)$
for any $v\in A$ if and only if
$$
|N_G(S)|\geq\sum\limits_{v\in S}{f(v)}-|S|+1
$$
for any nonempty subset $S\subseteq A$.

\section{The proof of Theorem 1.1}

\noindent{\it Proof of Theorem 1.1.} Suppose, to the contrary, that $G$ contains no spanning tree $T$ with $d_T(v)\geq k$ for any $v\in A$. By virtue
of Lemma 2.3, we conclude
\begin{align}\label{eq:3.1}
|N_G(S)|\leq(k-1)|S|
\end{align}
for some nonempty subset $S\subseteq A$. Choose a connected bipartite graph $G$ with partition $A\cup B$ such that its signless Laplacian spectral radius
is as large as possible, where $|A|=m$ and $|B|=n$. We claim that $S$ is a proper subset of $A$, that is, $S\subset A$. Otherwise, $S=A$. Then it follows
from \eqref{eq:3.1} and $S=A$ that $|N_G(A)|\leq(k-1)|A|=(k-1)m$. Combining this with $n\geq(k-1)m+1$, we deduce $n-|N_G(A)|\geq(k-1)m+1-(k-1)m=1$, which
is impossible because $G$ is connected. For convenience, we let $|S|=s$ and $|N_G(S)|=r$. Then we get $1\leq r\leq(k-1)s\leq(k-1)(m-1)$, and so
$1\leq s\leq m-1$. Obviously, there are no edges between $S$ and $B-N_G(S)$ in $G$. In terms of Lemma 2.1 and the choice of $G$ with bipartition $A\cup B$,
we conclude $G=K_{s,r}\nabla K_{m-s,n-r}$. It is clear that $G=K_{s,r}\nabla K_{m-s,n-r}$ is a spanning subgraph of $G_1=K_{s,(k-1)s}\nabla K_{m-s,n-(k-1)s}$.
Together with Lemma 2.1, we infer
\begin{align}\label{eq:3.2}
q(G)=q(K_{s,r}\nabla K_{m-s,n-r})\leq q(K_{s,(k-1)s}\nabla K_{m-s,n-(k-1)s}),
\end{align}
where the second equality holds if and only if $G=K_{s,(k-1)s}\nabla K_{m-s,n-(k-1)s}$. In what follows, we are to verify
$q(K_{s,(k-1)s}\nabla K_{m-s,n-(k-1)s})\leq q(K_{1,k-1}\nabla K_{m-1,n-k+1})$ with equality if and only if $s=1$.

If $s=1$, then $K_{s,(k-1)s}\nabla K_{m-s,n-(k-1)s}=K_{1,k-1}\nabla K_{m-1,n-k+1}$ and
$q(K_{s,(k-1)s}\nabla K_{m-s,n-(k-1)s})=q(K_{1,k-1}\nabla K_{m-1,n-k+1})$. Next, we consider $2\leq s\leq m-1$.

For the bipartite graph $G_1=K_{s,(k-1)s}\nabla K_{m-s,n-(k-1)s}$, the quotient matrix of the signless Laplacian matrix
$Q(G_1)=Q(K_{s,(k-1)s}\nabla K_{m-s,n-(k-1)s})$ by the partition $V(G_1)=S\cup(A-S)\cup N_{G_1}(S)\cup(B-N_{G_1}(S))$ is equal to
\begin{align*}
B_1=\left(
  \begin{array}{cccc}
    (k-1)s & 0 & (k-1)s & 0\\
    0 & n & (k-1)s & n-(k-1)s\\
    s & m-s & m & 0\\
    0 & m-s & 0 & m-s\\
  \end{array}
\right).
\end{align*}
Then the characteristic polynomial of $B_1$ is
\begin{align*}
\varphi_{B_1}(x)=&x^{4}-(2m+n+ks-2s)x^{3}+(m^{2}+mn+2kms+kns-3ms-ns-2ks^{2}+2s^{2})x^{2}\\
&+(kms^{2}-ms^{2}+kns^{2}-ns^{2}-km^{2}s+m^{2}s-kmns+mns)x.
\end{align*}
Notice that the partition $V(G_1)=S\cup(A-S)\cup N_{G_1}(S)\cup(B-N_{G_1}(S))$ is equitable. In view of Lemma 2.2, the largest root, say $q_1$, of
$\varphi_{B_1}(x)=0$ satisfies $q_1=q(G_1)=q(K_{s,(k-1)s}\nabla K_{m-s,n-(k-1)s})$. Note that $K_{m,(k-1)s}$ is a proper subgraph of
$G_1=K_{s,(k-1)s}\nabla K_{m-s,n-(k-1)s}$, and $G_1=K_{s,(k-1)s}\nabla K_{m-s,n-(k-1)s}$ is a proper subgraph of $K_{m,n}$. According to Lemma 2.1,
we have
\begin{align}\label{eq:3.3}
m+n=q(K_{m,n})>q_1=q(G_1)>q(K_{m,(k-1)s})=m+(k-1)s.
\end{align}

For the bipartite graph $G_*=K_{1,k-1}\nabla K_{m-1,n-k+1}$, the quotient matrix of $Q(G_*)$ in terms of the partition $V(G_*)=S\cup(A-S)\cup N_{G_*}(S)\cup(B-N_{G_*}(S))$ can be written as
\begin{align*}
B_*=\left(
  \begin{array}{cccc}
    k-1 & 0 & k-1 & 0\\
    0 & n & k-1 & n-k+1\\
    1 & m-1 & m & 0\\
    0 & m-1 & 0 & m-1\\
  \end{array}
\right),
\end{align*}
so its characteristic polynomial is
\begin{align*}
\varphi_{B_*}(x)=&x^{4}-(2m+n+k-2)x^{3}+(m^{2}+mn+2km+kn-3m-n-2k+2)x^{2}\\
&+(km-m+kn-n-km^{2}+m^{2}-kmn+mn)x.
\end{align*}
Note that the partition $V(G_*)=S\cup(A-S)\cup N_{G_*}(S)\cup(B-N_{G_*}(S))$ is equitable. According to Lemma 2.2, the largest root, say $q_*$,
of $\varphi_{B_*}(x)=0$ satisfies $q_*=q(G_*)=q(K_{1,k-1}\nabla K_{m-1,n-k+1})$.

Notice that $\varphi_{B_1}(q_1)=0$. By plugging the value $q_1$ into $x$ of $\varphi_{B_*}(x)-\varphi_{B_1}(x)$, we obtain
\begin{align}\label{eq:3.4}
\varphi_{B_*}(q_1)=\varphi_{B_*}(q_1)-\varphi_{B_1}(q_1)=q_1(s-1)\psi(q_1),
\end{align}
where $\psi(q_1)=(k-2)q_1^{2}-(2km+kn-3m-n-2ks-2k+2s+2)q_1-kms-km+ms+m-kns-kn+ns+n+km^{2}-m^{2}+kmn-mn$. In view of \eqref{eq:3.3} and $k\geq3$,
we easily see
\begin{align}\label{eq:3.5}
\psi(q_1)\leq\max\{\psi(m+n),\psi(m+(k-1)s)\}.
\end{align}

Recall that $k\geq3$ and $2\leq s\leq m-1$. We deduce
\begin{align}\label{eq:3.6}
\psi(m+n)=&(km+kn-m-n)s+km-m-n^{2}+kn-n-mn\nonumber\\
\leq&(km+kn-m-n)(m-1)+km-m-n^{2}+kn-n-mn\nonumber\\
=&-n^{2}+(km-2m)n+km^{2}-m^{2}.
\end{align}
Let $f(n)=-n^{2}+(km-2m)n+km^{2}-m^{2}$. Note that
$$
\frac{km-2m}{2}<(k-1)m<n
$$
by $k\geq3$ and $n\geq(k-1)m+1$. Then we deduce
\begin{align*}
f(n)<&f((k-1)m)\\
=&-(k-1)^{2}m^{2}+(km-2m)(k-1)m+km^{2}-m^{2}\\
=&0.
\end{align*}
Combining this with \eqref{eq:3.6}, we obtain
\begin{align}\label{eq:3.7}
\psi(m+n)\leq f(n)<f((k-1)m)=0.
\end{align}

By a direct computation, we get
\begin{align}\label{eq:3.8}
\psi(m+(k-1)s)=(k-1)(k(k-1)s^{2}-(kn-2k+2)s+m-n).
\end{align}
Let $h(s)=k(k-1)s^{2}-(kn-2k+2)s+m-n$. Recall that $2\leq s\leq m-1$. By a simple calculation, we have
\begin{align*}
h(2)=&4k(k-1)-2(kn-2k+2)+m-n\\
=&-(2k+1)n+4k^{2}+m-4\\
<&-(2k+1)(k-1)m+4k^{2}+m-4\\
=&-(2k^{2}-k-2)m+4k^{2}-4\\
\leq&-3(2k^{2}-k-2)+4k^{2}-4\\
=&-2k^{2}+3k-2\\
<&0
\end{align*}
and
\begin{align*}
h(m-1)=&k(k-1)(m-1)^{2}-(kn-2k+2)(m-1)+m-n\\
<&k(k-1)(m-1)^{2}-(k(k-1)m-2k+2)(m-1)+m-(k-1)m\\
=&-k(k-2)(m-1)-k+2\\
<&0
\end{align*}
due to $k\geq3$, $m\geq3$ and $n\geq(k-1)m+1>(k-1)m$. Thus, $h(s)\leq\max\{h(2),h(m-1)\}<0$ for $2\leq s\leq m-1$. Combining this with \eqref{eq:3.8},
we infer
\begin{align}\label{eq:3.9}
\psi(m+(k-1)s)=(k-1)h(s)<0.
\end{align}

Using \eqref{eq:3.5}, \eqref{eq:3.7} and \eqref{eq:3.9}, we conclude
\begin{align}\label{eq:3.10}
\psi(q_1)\leq\max\{\psi(m+n),\psi(m+(k-1)s)\}<0.
\end{align}

It follows from \eqref{eq:3.3}, \eqref{eq:3.4}, \eqref{eq:3.10}, $k\geq3$ and $2\leq s\leq m-1$ that
$$
\varphi_{B_*}(q_1)=q_1(s-1)\psi(q_1)<0,
$$
which yields $q(K_{s,(k-1)s}\nabla K_{m-s,n-(k-1)s})=q_1<q_*=q(K_{1,k-1}\nabla K_{m-1,n-k+1})$.

In conclusion, $q(K_{s,(k-1)s}\nabla K_{m-s,n-(k-1)s})\leq q(K_{1,k-1}\nabla K_{m-1,n-k+1})$ with equality if and only if $s=1$. Together with
\eqref{eq:3.2}, we obtain
$$
q(G)\leq q(K_{1,k-1}\nabla K_{m-1,n-k+1})
$$
with equality if and only if $G=K_{1,k-1}\nabla K_{m-1,n-k+1}$, which is a contradiction to the condition of Theorem 1.1 because
$G=K_{1,k-1}\nabla K_{m-1,n-k+1}$ has no spanning tree $T$ with $d_T(v)\geq k$ for any $v\in A$. This completes the proof of Theorem 1.1.   \hfill $\Box$

\medskip

\section*{Data availability statement}

My manuscript has no associated data.

\section*{Declaration of competing interest}

The authors declare that they have no known competing financial interests or personal relationships that could have appeared to influence the work
reported in this paper.

\section*{Acknowledgments}

This work was supported by the Natural Science Foundation of Jiangsu Province (Grant No. BK20241949). Project ZR2023MA078 supported by Shandong
Provincial Natural Science Foundation.

\end{document}